\documentclass[11pt,reqno]{amsart}

\usepackage{amssymb,amsmath}
\usepackage[dvips]{graphicx}
\usepackage{epsfig}
\usepackage{amssymb,amsmath}

\newcommand{\Diag}{\mbox{\rm Diag}}
\newcommand{\be}{\begin{eqnarray}}
\newcommand{\ee}{\end{eqnarray}}

\def\Ker{{\rm Ker}}

\def\wtil{\widetilde}

\newcommand{\clos}{\mbox{\rm clos}}

\newcommand{\es}{\emptyset}
\newcommand{\half}{\frac{1}{2}}

\newcommand{\eps}{{\varepsilon}}

\newcommand{\R}{{\mathbb R}}

\newcommand{\Nat}{{\mathbb N}}

\newcommand{\disk}{{\mathbb D}}

\newcommand{\Bk}{{\mathcal B}}
\newcommand{\Mk}{{\mathcal M}}

\newcommand{\bx}{{\bf x}}
\newcommand{\bob}{{\bf b}}
\def\bbe{{\bf e}}
\def\bd{{\bf d}}

\newcommand{\Nk}{{\mathcal N}}

\newcommand{\bo}{{\bf 0}}

\newcommand{\dist}{\mbox{\rm dist}}

\def\Om{{\Omega}}
\newcommand{\Lam}{{\Lambda}}
\def\Gam{{\Gamma}}
\newcommand{\lam}{\lambda}

\newcommand{\gam}{\gamma}

\newcommand{\Leb}{{\mathcal L}}

\renewcommand{\o}{\overline{1}}
\newcommand{\leb}{\mathcal{L}_2}

\def\ov{\overline}
\newtheorem{theorem}{Theorem}[section]
\def\Ok{\mathcal O}
\newtheorem{lemma}[theorem]{Lemma}
\newtheorem{cor}[theorem]{Corollary}
\newtheorem{prop}[theorem]{Proposition}

\theoremstyle{definition}

\newtheorem{example}[theorem]{Example}

\theoremstyle{definition}
\newtheorem{remark}[theorem]{Remark}

\numberwithin{equation}{section}

\input epsf.sty

\begin{document}

\thispagestyle{empty}

\title[zeros of power series]
{{ Zeros of $\{-1,0,1\}$ power series and \\[1.2ex]
connectedness loci for self-affine sets}}

\author{Pablo Shmerkin}

\address{Pablo Shmerkin,Box 354350, Department of Mathematics,
University of Washington, Seattle, WA 98195.
{\tt shmerkin@math.washington.edu}}

\author{Boris Solomyak}\thanks{
Both authors were supported in part by NSF grant
 \#DMS-0355187
}

\address{Boris Solomyak, Box 354350, Department of Mathematics,
University of Washington, Seattle, WA 98195.
{\tt solomyak@math.washington.edu}}

\begin{abstract}
We consider the set $\Om_2$ of double zeros in $(0,1)$ for power series
with coefficients in $\{-1,0,1\}$. We prove that $\Om_2$
is disconnected, and estimate $\min \Om_2$ with
high accuracy. We also show that $[2^{-1/2}-\eta,1)\subset \Om_2$
for some small, but explicit $\eta>0$ (this was only known for $\eta=0$).
These results have applications in the study of
infinite Bernoulli convolutions and connectedness properties of
self-affine fractals.
\end{abstract}
\date{\today}

\maketitle

\section{Introduction} \label{sec:intro}

Let
\begin{equation} \label{def-bk}
\Bk = \Bigl\{1+\sum_{n=1}^\infty a_n x^n:\ a_n \in \{-1,0,1\}\Bigr\}\,.
\end{equation}
We investigate the set
\begin{equation} \label{def-om2}
\Omega_2 = \{x\in (0,1):\ \exists\, f \in \Bk,\ f(x) = f'(x) = 0\}\,.
\end{equation}
Thus, $\Omega_2$ is the set of zeros of power series with coefficients
in $\{-1,0,1\}$ of order greater or equal to two (we call them
``double zeros'' for short). Since $\Bk$ is a normal family, $\Omega_2$
is relatively closed in $(0,1)$, so there exists $\alpha_2 =\min \Omega_2$.

We prove that the set $\Omega_2$ is disconnected and show that
$\alpha_2\approx .6684756$ (see Theorem~\ref{th-main}).
In fact, we found 58 distinct components in $\Omega_2$, and we conjecture
that there are infinitely many components. Numerical evidence indicates
that the structure of $\Om_2$ is very complicated.

Let $\wtil{\alpha}_2:= \sup((0,1)\setminus \Omega_2)$.
It is known that $[2^{-1/2},1) \subset \Omega_2$
(see Lemma \ref{lem-folk}), hence $\wtil{\alpha}_2\le 2^{-1/2}$.
We show that
$\wtil{\alpha}_2\le 2^{-1/2}-4\times  10^{-6}$ (see
Theorem~\ref{th-nonemptyint}). Numerical evidence suggests that
$\wtil{\alpha}_2 \approx 0.67$, but this is harder to prove rigorously.

The study of $\Omega_2$ is motivated by the work on infinite Bernoulli
convolutions \cite{Solerd,PeSo1,PeSo2,PeSchl} and on some fractal sets
\cite{JPoll}, where a key step is checking a certain {\em transversality
condition}. This condition holds precisely on $(0,1)\setminus \Omega_2$.
Those papers used the estimate $\alpha_2 \ge 0.649...$
obtained in \cite{Solerd} by computing the smallest double zero of a larger
class of power series $\wtil{\Bk}:=
\{1+\sum_{n=1}^\infty a_n x^n:\ a_n \in [-1,1]\}$.
Theorem~\ref{th-main} extends considerably (by more than 12\%)
the set of parameters where
the results of \cite{PeSo2,PeSchl,JPoll} apply.

Another motivation comes from the study of connectedness
loci for certain families of self-similar and self-affine fractals in the
plane. One of them is the ``Mandelbrot set for pairs of linear maps'' $\Mk$,
studied in \cite{BH,barn,bou,IJK,solcamb,bandt,sxu}.
We introduce two other connectedness loci, $\Nk$ and $\Ok$;
the latter one, associated with a linear map having a $2\times 2$ Jordan block,
coincides with $\Om_2 \cup (-\Om_2)$. Theorem~\ref{th-main} yields
new information about ``spikes,'' or ``antennas'' --- peculiar
features of $\Mk$ and $\Nk$.
Our second main result is Theorem~\ref{th-nonemptyint} in which
we obtain explicit neighborhoods of
(previously unknown) interior points of all three connectedness loci.

Let us make a few comments about the proofs. We use a C++ program, based on
a modification of Bandt's algorithm from \cite{bandt}, with
rigorous estimates, to rule out double zeros in specific intervals. This
program also indicates if there is a possible root in the interval, and
provides a polynomial which is the initial part of a power series in
$\Bk$ with a double zero in the interval. Once such a polynomial is found,
we use a simple argument (see Section \ref{sec:exist})
to prove the existence of the
function. It is completely rigorous, and its application
only uses {\sl Mathematica}
(or any similar package) to plot polynomials of degree up to $\approx 50$.
Thus, the lower estimate for $\alpha_2$ is computer-assisted in a more
substantial way than the upper estimate (which is just
``{\sl Mathematica}-assisted'').

In order to show that $\wtil{\alpha}_2<2^{-1/2}-4 \times 10^{-6}$ and obtain
neighborhoods of interior points in the connectedness loci, we use a
covering argument inspired by \cite{IJK} and \cite{sxu}. At one point we need
to check that a certain set is covered by the union of $3^5$ parallelograms,
which we do using a computer.

The paper is organized as follows. In Section 2 we provide the
background on iterated function systems and
discuss the relation between zeros of functions in $\Bk$ and
connectedness of self-affine fractals. We then state our results.
In Section 3 we show how to find double roots close to a local minimum of
a polynomial with certain properties.
In Section 4 we establish the covering results and estimate $\wtil{\alpha}_2$.
In Section 5 we prove the existence of gaps in $\Omega_2$ and estimate
$\alpha_2$.
Section 6 is devoted to some variants and generalizations.
Section 7 contains proofs of several auxiliary results.


\section{Preliminaries on IFS and statement of results}

An {\em iterated function system} (IFS) is a finite collection of
(strict) contractions $\{f_1,\ldots,f_m\}$ on a complete metric space. Given
such a system, there is a unique nonempty compact set $E$ satisfying
$E = \bigcup_{i\le m} f_i(E)$, called the {\em attractor} of the IFS,
see \cite{hutch}.
We only consider IFS on $\R^d$ of the form $\{T_i \bx + \bob_i\}_{i\le m}$
where $T_i$ are linear maps and $\bob_i  \in \R^d$. Their attractors are
called {\em self-affine}. For the maps to be contractive (in some
norm) it is necessary and sufficient that all the eigenvalues of $T_i$ are
less than 1 in absolute value.

We investigate when attractors are connected in the simplest case
$m=2$. The following is well-known.

\begin{prop}[see \cite{Hata}] \label{prop-hata}
The attractor $E$ of an IFS $\{f_1,f_2\}$ is connected if and only if
$f_1(E) \cap f_2(E) \ne \es$.
{\em (Of course, the ``only if'' direction is obvious.)}

\end{prop}

For IFS of two
affine maps there is a simple sufficient condition for connectedness.
We can assume $\bob_1=\bo$ without loss of generality, making
a change of variable.

\begin{lemma}[folklore] \label{lem-folk} Let $\{T_1\bx,  \,T_2\bx + \bob\}$
be an IFS of contracting affine maps, such that $\max\{\|T_1\|,\|T_2\|\}<1$ in
some operator norm, and
$|\det(T_1)| + |\det(T_2)| \ge 1$. Then the attractor is connected.
\end{lemma}

We include a proof
for completeness (see Section 7).
Next we specialize even more, assuming that $T = T_1 = T_2$, and state
a criterion for connectedness in terms of zeros of power series.

Let $E=E(T,\bob)$ be the attractor of the IFS $\{T\bx, T\bx + \bob\}$,
{\em i.e.}, the unique nonempty compact set in $\R^d$ satisfying
\be \label{saf1}
E = TE \cup (TE + \bob).
\ee
Observe that
\be \label{serep}
E(T,\bob) = \Bigl\{\sum_{n=0}^\infty a_n T^n \bob:\ a_n \in \{0,1\} \Bigr\}
\ee
since the right-hand side is well-defined and satisfies
(\ref{saf1}).

We can assume, without loss of generality, that
$\bob$ is a cyclic vector
for $T$, that is, $H  := Span\{T^k \bob:\ k\ge 0\} = \R^d$.
Indeed, otherwise
we can replace $T$ by the restriction of $T$ to $H$
and consider the corresponding IFS on $H$.

Combining Proposition~\ref{prop-hata}
and (\ref{serep}) easily implies the following criterion,
which is known, at least in special cases.
Recall that $\Bk$ is defined in (\ref{def-bk}).

\begin{prop} \label{prop-connect} Let $T$ be a linear contraction with
(possibly complex) eigenvalues $\lam_j$, for $j=1,\ldots, m$,
having algebraic
multipicities $k_j \ge 1$, and geometric multiplicities equal
to one. Let $\bob$ be a cyclic vector for $T$. Then $E(T,\bob)$ is
connected if and only if there exists $f \in \Bk$ such that
\be \label{eq-zero}
f(\lam_j) = \ldots = f^{(k_j-1)}(\lam_j) = 0,\ \ \ j=1,\ldots,m.
\ee
In particular, connectedness does not depend on $\bob$.
\end{prop}

Since this is a key statement relating connectnedness
of self-affine
sets to zeros of power series, we include a short proof in Section 7.
Combining Proposition~\ref{prop-connect} and Lemma~\ref{lem-folk}
yields the following result.
We denote by $\disk$ the open unit disk in the complex plane.

\begin{cor} \label{cor-connect1}
Let $\Lam = \{\lam_1,\ldots,\lam_m\}\subset \disk$ and let $k_j = k(\lam_j)
\ge 1$ be such that for any nonreal $\lam\in \Lam$, we have
$\ov{\lam} \in \Lam$ and $k(\ov{\lam}) = k(\lam)$.
If $\prod_{j=1}^m |\lam_j|^{k_j} \ge 1/2$,
then there exists $f\in \Bk$ having zeros at $\lam_j$ of multiplicity
$\ge k(\lam_j)$ for $j=1,\ldots,m$.
\end{cor}

In particular, we obtain that for any $k\ge 1$, every $\lam \in [2^{-1/k},1)$
is a zero of multiplicity $\ge k$ for some power series
in $\Bk$. In \cite[Section 3]{BBBP} it is
asked whether there exist power series
(or polynomials) with coefficients in $\{-1,0,1\}$
having a $k$-th order root strictly inside the unit circle for arbitrary $k$.
Corollary~\ref{cor-connect1} answers the question for power series in a strong
quantitative way, but the question
for polynomials is much harder and remains open.

\medskip

From now on, we restrict ourselves to the case $d=2$.
Applying an
invertible linear transformation as a conjugacy, we can assume without
loss of generality that $T$ is one of the following:
$$
{\rm (i)}\ \ T = \left(\begin{array}{rr} a & b \\
-b & a \end{array} \right),
\ \ \ \ \ {\rm (ii)}\ \ T = \left(\begin{array}{cc} \gam & 0 \\ 0 & \lam
\end{array} \right),
\ \ \ \ \ {\rm (iii)}\ \ T = \left(\begin{array}{cc} \lam & 1 \\ 0 & \lam
\end{array} \right),
$$
where $a, b, \lam, \gam$ are real,
$a^2 + b^2 < 1$, $b\ne 0$, $|\lam|, |\gam| < 1$, and  $\gam \ne \lam$.
Each of the cases leads to a set which we call the {\em connectedness
locus} for the corresponding family of self-affine sets. Namely, we consider
the sets
\begin{eqnarray*}
\Mk & := & \{z=a+ib\in \disk:\ \exists\,f\in \Bk,\ f(z)=0\},\\
\Nk & := & \{(\gam,\lam) \in (-1,1)^2:\ \exists\,f\in \Bk,\ f(\gam) = f(\lam)
=0\},\\
\Ok & := & \{\lam\in (-1,1):\ \exists\,f\in \Bk,\ f(\lam) = f'(\lam)=0\}.
\end{eqnarray*}
Thus, $\Mk,\Nk,\Ok$ are essentially the sets
of parameters for which the attractors
in cases (i),(ii),(iii) respectively are connected.
(It is natural to allow $\gam=\lam$
in $\Nk$ and $b=0$ in $\Mk$ to ensure that the sets are relatively closed in
$\disk$.)
By Lemma~\ref{lem-folk},
\begin{eqnarray*}
\Mk & \supset & \Mk_t:= \{ \lambda\in\mathbb{D}:\
 |\lambda| \ge 2^{-1/2}\},\\
\Nk & \supset & \Nk_t: = \{ (\gamma,\lambda)\in(-1,1)^2:\
 |\gamma\lambda| \ge 1/2\},\\
\Ok & \supset & \Ok_t: = \{ \lambda\in(-1,1):\
|\lambda|\ge 2^{-1/2}  \}.
\end{eqnarray*}
We refer to $\Mk_t, \Nk_t, \Ok_t$ as ``trivial parts''
of the corresponding sets.

The set $\Mk$
was studied by several authors, see
\cite{BH,barn,bou,solcamb,bandt,sxu}. In particular, Bousch \cite{bou}
proved that $\Mk$ is connected and locally connected. The set $\Nk$
has not been studied as much, although partial results are obtained in
\cite{SolSAF}, where it is shown that a large ``chunk'' of $\Nk$ is
connected (all of $\Nk$ is conjectured to be connected).
An approximation to
$\Nk\cap (0,1)^2$ is depicted in Figure 1, with the non-trivial
part shown in black. The picture is created with a program of C. Bandt;
note that the visible
disconnected pieces are a computing artefact.


\begin{figure}[ht]
\scalebox{.5}{\epsfig{figure=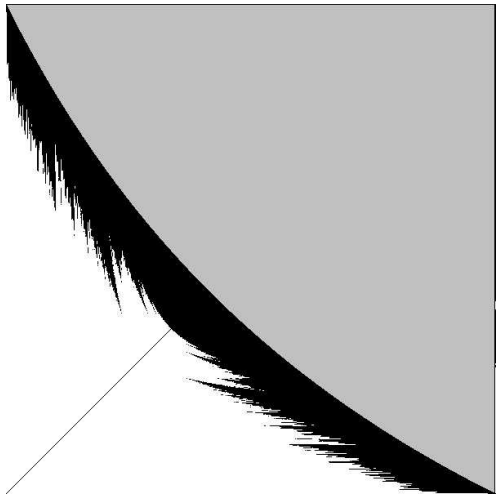,height=13cm}}
\caption{Part of the connectedness locus $\Nk$}
\end{figure}

By symmetry, we have $\Ok = \Om_2 \cap (-\Om_2)$. The set $\Om_2$,
defined in (\ref{def-om2}), is our main object of study.
Since $\Bk$ is a normal family, $\Omega_2$
is relatively closed in $(0,1)$, so there exists $\alpha_2 =\min \Omega_2$.

\begin{theorem} \label{th-main}
{\bf (i)} $\alpha_2  \in (0.6684755,0.6684757)$;

{\bf (ii)} $\Omega_2$ is disconnected. In fact, the intervals $I_j$, for
$j\le 5$, lie in
distinct components of $(0,1)\setminus \Omega_2$ where
$I_j = .668+10^{-3} I_j'$, and
$I_1' = (.478,.489), \ I_2' = (.632,.654),\ I_3' = (1.282,1.306),\
I_4' = (1.327,1.333),\ I_5' = (1.343,1.352)$.
\end{theorem}

This theorem is proved in Section 5, using results of
Section 3 for part of the proof.

\begin{remark} Numerical evidence suggests that there are infinitely
many components of $(0,1)\setminus \Omega_2$.
We do not have a proof of that.
The topological structure of $\Omega_2$ appears to be very complicated.
It can be proved rigorously that the five ``gaps'' above are the largest.
\end{remark}

\begin{remark}
Let $\wtil{\alpha}_2:= \sup((0,1)\setminus \Omega_2)$.
We have $[2^{-1/2},1) \subset \Omega_2$ by Lemma \ref{lem-folk},
hence $\wtil{\alpha}_2\le 2^{-1/2}$. We are able to show in fact that
$\wtil{\alpha}_2\le 2^{-1/2}-4\times  10^{-6}$ (see
Theorem~\ref{th-nonemptyint}
below). Our rigorous
numerical results also yield $\wtil{\alpha}_2\ge .669355$, and it seems that
$\wtil{\alpha}_2\le .67$, but we do not have a proof.
\end{remark}

Returning to the sets $\Mk$, $\Nk$, we note that they are related to
$\Ok$. In fact, we have the following, denoting for $F\subset \R$,\
$\Diag(F): = \{(\lam,\lam):\ \lam\in F\}$.

\begin{lemma} \label{lem-anten} {\bf (i)} $\clos(\Mk\setminus \R) \cap
\R \subset \Ok$;

{\bf (ii)} $\clos(\Nk\setminus \Diag(\R)) \cap \Diag(\R) \subset \Diag(\Ok)$.
\end{lemma}

This follows by an easy compactness argument.
See \cite[Lemma 2.5]{solcamb} for the proof of
(i); (ii) is proved in Section 7.

As a consequence of Lemma~\ref{lem-anten},
the pictures of $\Mk$ near the real axis reveal something
about the structure of $\Ok$. In fact, Figure 7
in \cite{bandt} served as an inspiration for our work.
Lemma~\ref{lem-anten} also explains an interesting feature of $\Mk$ and
$\Nk$, namely the ``antennas.'' The antenna for $\Mk$ (let us denote it
by $\Gam(\Mk)$) is defined as the connected component of $[\half,1) \setminus
\clos(\Mk\setminus \R)$ containing $\half$
(there is obviously a symmetric antenna on
the negative real axis); it was first discovered in \cite{BH} and studied
in \cite{solcamb,bandt}. Similarly, the antenna for $\Nk$, denoted by
$\Gam(\Nk)$, may be defined as the connected component of
$\Diag([\half,1))\setminus \clos(\Nk\setminus \Diag(\R))$ containing $(\half,
\half)$.
Bandt \cite{bandt} noted that the ``tip of the antenna''
$\sup\Gam(\Mk)$ is the infimum of the set of
double zeros of
power series in $\Bk$ with infinitely many coefficients not equal to $+1$.
It is not hard to show that
$\sup\Gam(\Nk)$ is the infimum of the set of
double zeros of
power series in $\Bk$ with infinitely many coefficients not equal to $-1$.
It seems very likely that $\sup\Gam(\Mk)=\sup\Gam(\Nk)=\alpha_2$, but we
do not know how to prove this. However,
as a by-product of our investigation, we obtain the following corollary
proved in Section 5.

\begin{cor} \label{cor-ant}
We have  $\sup\Gam(\Mk),\ \sup\Gam(\Nk)\in (0.6684755,0.6684757)$.
\end{cor}

Our second main result concerns non-trivial interior points of the
connectedness loci.
In \cite{IJK} and \cite{sxu} some chunks of interior points of
$\Mk \backslash \Mk_t$ were found.
Although numerical experimentation indicates that $\Nk\backslash
\Nk_t$ and $\Ok\backslash\Ok_t$
also have nonempty interior, this had not been proved rigorously before.

\begin{theorem} \label{th-nonemptyint}
Let $\eta=4\times 10^{-6}$. Then
\begin{eqnarray}
\Mk & \supset & \Delta_1 =  \Bigl\{ (\gamma,\lambda)\in(0,1)^2:\
 |\gamma-2^{-1/2}|,|\lambda-2^{-1/2}|< \eta\Bigr\},\label{eq-mandelint1}\\
\Nk & \supset & \Delta_2 = \Bigl\{ \lambda\in\mathbb{D}:\
|\lambda-2^{-1/2}|< 3\eta/4 \Bigr\},\label{eq-mandelint2}\\
\Ok & \supset & \Delta_3 = \Bigl\{ \lambda\in(0,1):\
 |\lambda-2^{-1/2}|<\eta  \Bigr\}\label{eq-mandelint3}.
\end{eqnarray}
\end{theorem}

This theorem is proved in Section 5.
The value of $\eta$ is very small and clearly not optimal, but nevertheless
it is an explicit constant.


\section{Existence of double roots I} \label{sec:exist}

We denote by $\Bk_n$ the
subset of $\Bk$ consisting of polynomials
having degree less than or equal to $n$.

Let us say that $(P,n,a,b)$ is {\em good} if $P \in \Bk_n$, $0.5 < a < b < 1$
(in reality, we will only consider $0.66 < a < b < 0.68$),
\begin{equation} \label{eq1}
P(a) > a^{n+1}/(1-a),\ \ \ P(b) > b^{n+1}/(1-b),
\end{equation}
$P(x) > 0 $ for all $x\in [a,b]$, and
\begin{equation} \label{eq2}
\exists\,x\in (a,b):\ P(x) < x^{n+1}/(1-x).
\end{equation}

\begin{lemma} \label{lem-good} Suppose that $(P,n,a,b)$ is good.
Let $Q(x) = P(x) - x^m$ where $m$ is the minimal integer greater or equal to
$n+1$, such that
$Q(x) > 0$ on $[a,b]$. Then $(Q,m,a,b)$ is good.
\end{lemma}

{\em Proof.} It is clear that $Q \in
\Bk_m$. We easily check that
$$
Q(a) = P(a) - a^m > a^n/(1-a) - a^m > a^{m+1}/(1-a),
$$
since $a \in (0,1)$, and similarly, $Q(b) > b^{m+1}/(1-b)$. It remains to
check the last condition.
Either $m=n+1$, in which case we note that
$$
Q(x) = P(x) -x^{n+1} < x^{n+1}/(1-x) - x^{n+1} = x^{n+2}/(1-x),
$$
or $m > n+1$, in which case there exists $t\in (a,b)$ such that
$P(t) - t^{m-1} \le 0$.  Then
$$
Q(t) = P(t) - t^m = (P(t) - t^{m-1}) + (t^{m-1}-t^m) \le t^{m-1}-t^m <
t^{m+1}/(1-t),
$$
since $t$ is greater than $\half$. Clearly $t\ne a,b$,
and the proof is complete. \qed

\begin{cor} \label{cor-exist}
Suppose that $(P,n,a,b)$ is good.
Then there exists $f\in \Bk$ such that $P$ is the initial part of $f$, and
$f$ has a double zero in $(a,b)$.
\end{cor}

{\em Proof.} Iterating Lemma~\ref{lem-good} we obtain a sequence of
polynomials $Q_j \in \Bk_{m_j}$ such that $m_j\to \infty$ and
$Q_j$ is the initial part of $Q_{j+1}$ for all $j$. Then $Q_j \to f\in \Bk$
uniformly on compact subsets of the unit disk in the complex plane.
Since $Q_j$ are all positive on $[a,b]$, we have that $f(x) \ge 0$ on $[a,b]$.
Since $Q_j(x_j) < x_j^{m_j}$ for some $x_j \in (a,b)$, we have that
$\min_{[a,b]} f = 0$. On the other hand, (\ref{eq1}) implies that
any function in $\Bk$ with initial part $P$ is strictly positive at $a$ and $b$.
It follows that $f$ has a zero in $(a,b)$
of order at least two.
\qed

\begin{remark} \label{rem-cor}
The function $f(x)=1+\sum_{j=1}^\infty a_j x^j$ obtained in
Corollary~\ref{cor-exist} has the property that $a_j \in \{0,-1\}$ for all
$j \ge n+1$ and there are infinitely many 0's and $-1$'s.
The only thing to check is that there are infinitely many 0's; the rest
is obvious by construction in Lemma~\ref{lem-good}. Suppose that
$a_j = -1$ for all $j \ge N$ for some $N\in \Nat$. Then at some step in
our construction we have good $(Q,N,a,b)$, which implies $Q(x) < x^{N+1}/(1-x)$
for some $x\in (a,b)$. Then $f(x) = Q(x) - x^{N+1}/(1-x) < 0$, which is
a contradiction.
\end{remark}

Corollary~\ref{cor-exist} is not very efficient numerically, since it allows us
to find a double zero with an error of order $b^{n/2}$.
The next statement shows that this can be improved considerably.

\begin{cor} \label{cor-exist2}
Assume that $(P,n,a,b)$ is good, $b\le .68$, and $n>10$. Further,
suppose that $P'$ has a zero at $y\in (a,b)$, and
\begin{equation} \label{eq-sder}
P''(x) > C''\ \ \ \mbox{for all}\ x\in [a,b]
\end{equation}
Then there exists $f\in\Bk$ such that $P$ is the initial part of $f$ and
$f$ has a double root in the interval $(y-\eta,y+\eta)$, where
\[
\eta = \frac{1+(1-b) (n+1)}{C'' (1-b)^2}\, b^{n+1} < (4/C'') (n+1) b^{n+1}.
\]
\end{cor}
{\em Proof.} From Corollary \ref{cor-exist} we know that
there exists $f\in\Bk$ such
that $P$ is its initial part, and $f$ has a double root in $(a,b)$;
let $r$ be this root. Then $P'(y)=f'(r)=0$
whence, using the intermediate value theorem,
\begin{equation} \label{eq-lb}
| P'(r) - f'(r) | = | P'(r) - P'(y) | > C'' |r-y|.
\end{equation}
Note however that
\begin{equation} \label{eq-ub}
| P'(r) - f'(r) | \le \sum_{i=n+1}^\infty i r^{i-1} =
\frac{(n+1)r^n(1-r) + r^{n+1}}{(1-r)^2} <
\frac{1+(1-b)(n+1)}{(1-b)^2}\, b^n.
\end{equation}
Combining (\ref{eq-lb}) and (\ref{eq-ub}) yields the corollary. \qed

\begin{example} \label{ex-double}
Consider
\begin{eqnarray*}
P(x) & = &
1 - x^{1} - x^{2} - x^{3} + x^{4} + x^{6} + x^{7} + x^{9} + x^{10} + x^{12}+x^{13} \\& + &
   x^{14} + x^{16} + x^{17} + x^{18} + x^{20} + x^{21} + x^{23} + x^{24} + x^{25} \\& +&
   x^{26} + x^{27} + x^{28} + x^{29} + x^{3 1} + x^{32}+ x^{33} + x^{34} + x^{35}\\& +&
   x^{36} + x^{37} - x^{38} + x^{39} + x^{40}+x^{41} + x^{42} + x^{43} - x^{44}\\ & + &
   x^{45} - x^{46} - x^{47} + x^{48} - x^{49} - x^{50}\in\mathcal{B}_{50}.
\end{eqnarray*}
Let $a=0.668470, b= 0.668482$.
We claim that
$(P,50,a,b)$ is good. We checked that
$P(x) > 0$ on $[a,b]$ by plotting the graph of $P$ on $[a,b]$ in
{\sl Mathematica}; see Figure \ref{fig-good}.

\begin{figure}
\centering
\includegraphics[width=0.9\textwidth]{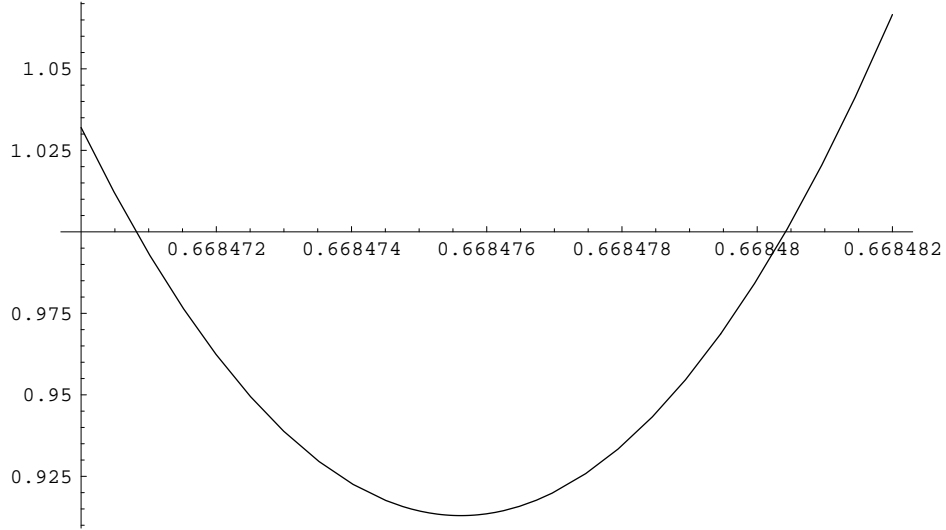}
\caption{Checking the existence of double roots: the polynomial $P$ was obtained with the help of a C++ program. We plotted $P(x)(1-x)/x^{51}$ on the interval $(0.668470,0.668482)$; it is clear from the picture that $(P,51,0.668470,0.668482)$ is good.}
\label{fig-good}
\end{figure}

We have
$$
P(a)(1-a)a^{-51} \approx 1.03199,\ \ \ \ P(b)(1-b) b^{-51} \approx 1.06665,
$$
verifying (\ref{eq1}). On the other hand,
$$
P(r)(1-r) r^{-51} \approx 0.912958, \ \ \ \mbox{for}\ r=0.6684756,
$$
\begin{sloppypar}
verifying (\ref{eq2}). Thus, Corollary~\ref{cor-exist} applies,
and we have a double zero in $[0.668470,0.668482]$.
\end{sloppypar}

Using Corollary~\ref{cor-exist2}, we obtain a more precise estimate.
We have
\begin{eqnarray}
P(r+2\times 10^{-8}) - P(r) & \approx & 8.89847\times 10^{-15}>0 \nonumber,\\
P(r-2\times 10^{-8}) - P(r) & \approx & 2.04733\times 10^{-15}>0 \nonumber,
\end{eqnarray}
which implies that there exists $y\in (0.66847558,0.66847562)$ such that
$P'(y)=0$. We checked that $P''(x) > 20$ on $[a,b]$
by plotting the graph of $P''$ on $[a,b]$ in
{\sl Mathematica}.
Thus, Corollary~\ref{cor-exist2} applies. Since
$50 b^{50}\cdot (4/20) < 2\times 10^{-8}$, we obtain that there exists
$f\in \Bk$ with a double zero in $(0.66847556,0.66847564)$.
\end{example}


\section{Existence of double roots II and connectedness
loci} \label{sec:exist2}

Here we prove Theorem \ref{th-nonemptyint}.
The proof is based on several lemmas.

\begin{lemma} \label{lemma-attrconn1}
Let $T$ be a contracting linear map on $\mathbb{R}^2$,
and let $\bob\in\mathbb{R}^2$ be a vector such that
$T^k\bob, k\ge 0$, span $\mathbb{R}^2$.
Denote the attractor of $\{T\bx, T\bx + \bob\}$ by $E$.
Let $\mathbf{v}$ be a point of the form
\[
\mathbf{v} = \sum_{i=1}^k a_i T^{-i} \bob,
\]
where $a_i\in\{-1,0,1\}$ and $a_k=1$. Assume that there exists a set $U$ containing $\mathbf{v}$ and $j\in\mathbb{N}$ such that
\begin{equation} \label{eq-largerifs}
U \subset \bigcup_{u\in\{-1,0,1\}^j}
(T+u_1 \bob)\circ\cdots\circ (T+u_j \bob)(U).
\end{equation}
Then $K$ is connected.
\end{lemma}
{\em Proof}. There are similar results in e.g.\ \cite{IJK}, but we sketch a
proof for completeness. Let $\wtil{E}$ denote the attractor of the IFS
$\{T\bx-\bob, \,T\bx,\, T\bx + \bob\}$. Condition (\ref{eq-largerifs})
implies that $U$ is contained in $\wtil{E}$. In particular,
$\mathbf{v}\in\wtil{E}$. Recalling the form of $\mathbf{v}$ we get that
\[
\sum_{i=1}^k a_i T^{-i} \bob = \sum_{i=0}^\infty c_i T^i \bob,\quad
\mbox{for some}\ c_i\in\{-1,0,1\}.
\]
Applying $T^k$ from the left on both sides we obtain a power series
$f\in\mathcal{B}$ such that $f(T)\bob=0$. Write $f=f_1-f_2$, where $f_i$ has
coefficients $0, 1$ only. In particular, since $a_k\ne 0$, the power series
$f_1$ and $f_2$
have different constant terms, hence
\[
f_1(T)\bob = f_2(T)\bob \in T E \cap (T E + \bob).
\]
By Proposition~\ref{prop-hata}, the set
$E$ is connected.  \qed

\begin{lemma} \label{lemma-parallelcover}
Let
\[
T = \left(
\begin{array}{cc}
  2^{-1/2} & 0.7 \\
  0 & 2^{-1/2} \\
 \end{array}
\right);\quad \mathbf{b}=\left(%
\begin{array}{c}
  1 \\
  1 \\
 \end{array}%
\right).
\]
Let also $\mathbf{p}=(-2.10,0.20), \mathbf{q}=(4.90,2.45)$.
Denote by $U$ the open parallelogram with the
vertices $\pm \mathbf{p}, \pm \mathbf{q}$, and denote by $V$ the open
parallelogram with the vertices $\pm 0.95 \mathbf{p},  \pm 0.95\mathbf{q}$.
Then
\[
\mathbf{v} = T^{-1}\bob - T^{-2}\bob - T^{-3}\bob - T^{-4}\bob + T^{-5}\bob
\approx (2.95837, 1.75736) \in V,
\]
and
\begin{equation} \label{eq-parallcover}
U \subset \bigcup_{u\in\{-1,0,1\}^5} (T+u_1\bob)\cdots(T+u_5\bob)(V).
\end{equation}
\end{lemma}

\begin{figure}
\centering
\includegraphics[width=0.9\textwidth]{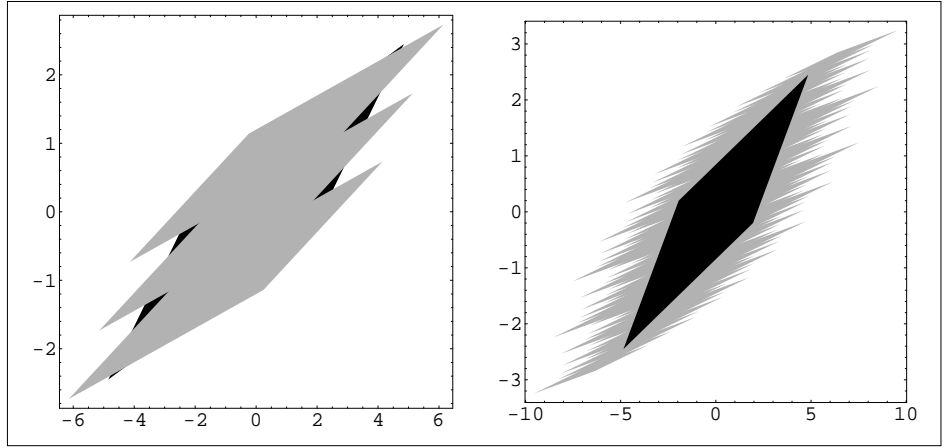}
\caption{The parallelogram $U$ is pictured in black. On the left, the gray figure is $(T V-\bob)\cup (T V) \cup (T V+\bob)$; notice that $U$ comes close to being covered by the iterates of $V$ already in the first step. However, $5$ iterations steps are needed in order to get a complete covering. This is illustrated by the picture on the right; there, the gray figure corresponds to the right-hand side in (\ref{eq-parallcover}).}
\label{fig-covering}
\end{figure}

{\em Proof}. This is the part of the proof that is computer assisted. The parallelograms $U$ and $V$ and the vector $\mathbf{v}$ were obtained through experimentation with {\em Mathematica}.

The coordinates of $\mathbf{v}$ in the base $\{\mathbf{p},\mathbf{q}\}$ are $(\alpha,\beta) \approx (0.223,0.699)$. Since $|\alpha|+|\beta|<0.95$, we get that $\mathbf{v}\in V$. We checked (\ref{eq-parallcover}) rigorously using an algorithm which we now describe.

Given a parallelogram $P$, a collection of parallelograms $\mathcal{Q}$ and a depth $n$, we wrote a routine with the following pseudo-code:

\begin{itemize}
\item[(i)] If $P$ is completely covered by some $Q\in\mathcal{Q}$, return {\em true} and exit (by convexity, it is enough to check whether all $4$ vertices of $P$ are contained in some $Q$).
\item[(ii)] If $n=0$ return {\em false} and exit.
\item[(iii)]Split $P$ into $4$ congruent pieces $P_i$. For each $i$, compute $\mathcal{Q}_i$, the family of parallelograms in $\mathcal{Q}$ intersecting $P_i$.
\item[(iv)] If any $\mathcal{Q}_i$ is empty, return {\em false} and exit.
\item[(v)] For each $i$, run the routine with input $P_i$, $\mathcal{Q}_i$ and $n-1$. If all of them return {\em true}, then return {\em true}; otherwise, return {\em false}.
\end{itemize}

Thus if the routine returns {\em true} then $P$ is covered by the union of the family $\mathcal{Q}$. This was indeed the case with the input $P=U$, $n=8$ and
\[
\mathcal{Q}= \left\{  (T+u_1\bob)\cdots(T+u_5\bob)(V) : u\in\{-1,0,1\}^5 \right\},
\]
and this proves the lemma; Figure \ref{fig-covering} depicts the situation
graphically. \qed

\medskip

{\em Proof of Theorem \ref{th-nonemptyint}}. Let
$T,\bob, \mathbf{p},\mathbf{q}, U, V$ be as in Lemma \ref{lemma-parallelcover}.
Using Lemma \ref{lemma-parallelcover} we see that for sufficiently small $\eta$,
\begin{equation} \label{eq-claimcovering}
\|T'-T\|<\eta \quad \Longrightarrow \quad V \subset
(T'+u_1 \bob)\circ\cdots\circ (T'+u_5 \bob)(V).
\end{equation}
The bulk of the proof will consist in obtaining an explicit value of $\eta$
such that (\ref{eq-claimcovering}) holds.

\begin{lemma} \label{lem-vsp1} The condition (\ref{eq-claimcovering}) holds
for $\eta = 4\times 10^{-6}$.
\end{lemma}

The proof of the lemma is straightforward, but technical, so we postpone it till
Section 7. As mentioned above, the existence of such positive $\eta$ is
obvious.

An immediate consequence of (\ref{eq-claimcovering}) and Lemmas
\ref{lemma-attrconn1} and \ref{lemma-parallelcover} is
\begin{equation} \label{eq-claimconnected}
\|T'-T\|<\eta \quad \Longrightarrow \textrm{ the attractor of }
\{T'\bx,T'\bx+\bob\} \textrm{ is connected}.
\end{equation}

It is convenient to use the $\ell^\infty$ norm in $\R^2$ and
the operator matrix norm; recall that the latter is computed as
the maximum of $\ell^1$ norms of the rows. Let
\begin{eqnarray}
M_1(\gamma,\lambda) & = & \left(%
\begin{array}{cc}
  \gamma & 0.7 \\
  0 & \lambda \\
\end{array}%
\right),\quad \gamma\ne\lambda,\nonumber\\
M_2(r,\varepsilon) & = &  \left(%
\begin{array}{cc}
  \rho & 0.7 \\
  -\varepsilon & \rho \\
\end{array}%
\right), \quad \rho,\varepsilon>0, \nonumber\\
M_3(\lambda) & = &  \left(%
\begin{array}{cc}
  \lambda & 0.7 \\
  0 & \lambda \nonumber
\end{array}%
\right).
\end{eqnarray}
Note that $M_1(\gamma,\lambda)$ is conjugate to the diagonal matrix with
eigenvalues $\gamma$ and $\lambda$; likewise, $M_3(\lambda)$ is conjugate to
the standard Jordan block with eigenvalue $\lambda$. Note also that
\begin{eqnarray}
\|M_1(\gamma,\lambda) - T \| & = & \max(|2^{-1/2}-\gamma|,|2^{-1/2}-\lambda|);
\nonumber \\
\|M_3(\gamma,\lambda) - T \| & = & |2^{-1/2}-\lambda|\nonumber.
\end{eqnarray}
From this, Lemma~\ref{lem-vsp1} and (\ref{eq-claimconnected})
we obtain (\ref{eq-mandelint1})
and (\ref{eq-mandelint3}).

It remains to consider the complex eigenvalue case. The matrix
$M_2(\rho,\varepsilon)$ has eigenvalues $\rho \pm i\sqrt{0.7 \varepsilon} $.
Let $\lambda$ be a non-real complex number such that $|\lambda-2^{-1/2}|< (3/4)
 \eta$; without loss of generality assume that $\textrm{Im}(\lambda)>0$, and
write $\lambda = \rho + i \sqrt{0.7 \varepsilon}$. We have that
$|\rho-2^{-1/2}| < (3/4) \eta$ and
\[
\sqrt{0.7 \varepsilon} < (3/4) \eta \,\ \Rightarrow\ \, \varepsilon <
\frac{9}{16 \times 0.7} \,\eta^2 < (1/4)\eta.
\]
Therefore $|\rho-2^{-1/2}|+\varepsilon < \eta$, and from this we conclude that
$\|M_2(\rho,\varepsilon)-T\|<\eta$. Invoking (\ref{eq-claimconnected}) again,
this shows that (\ref{eq-mandelint2}) is verified, which completes the proof.
\qed

\begin{remark}
Although the proof of Theorem \ref{th-nonemptyint} is inspired by analogous
results for self-similar sets which appeared in  \cite{IJK} and \cite{sxu},
the more complicated geometry of self-affine sets introduce some additional
difficulties. For example, the use of a computer to first find the sets $U, V$
and the vector $\mathbf{v}$,
and then for checking that the covering property holds,
becomes essential (in the self-similar case, some attractors are actually
rectangles, which allows to do a purely algebraic analysis in certain region,
see \cite{sxu}).  The chunks of interior points of $\Mk\setminus \Mk_t$
which were found in those papers are away from the real line.
\end{remark}

\begin{remark}
Let $K_{\gamma,\lambda}$ be the attractor of $\{T_{\gamma,\lambda},\,
T_{\gamma,\lambda}+(1,1)\}$, where $T_{\gamma,\lambda}$ is a diagonal map with
eigenvalues $\gamma,\lambda$. In \cite{shm} the almost sure Hausdorff dimension
 of $K_{\gamma,\lambda}$ was found in the region $(0,1)^2\backslash\Nk_t$.
The result was new only in the region $\Nk\backslash\Nk_t$; hence
Theorem \ref{th-nonemptyint} makes that result effective, by showing that
$\leb(\Nk\backslash\Nk_t)>0$.
\end{remark}

\begin{remark}
The value of $\eta$ found in Lemma~\ref{lem-vsp1} is extremely small, but
graphical experimentation suggests that the same covering argument, even with
the same covering, works for a large range of parameters. We believe that it
should be possible to extend the result to show that actually $(0.7,2^{-1/2})
\subset \Om_2$ (we recall that computer results suggest that indeed
$(0.67,2^{-1/2}) \subset \Om_2$, but this appears harder to prove rigorously).
\end{remark}


\section{{\bf Existence of gaps}} \label{sec:gaps}

We describe an algorithm which we use to rigorously prove the existence of gaps in the set $\Omega_2$. As a by-product, the results of Section \ref{sec:exist} yield bounds on how large the gaps may be. In practice, we are able to obtain an accurate description of the set $\Omega_2$ up to an error of $10^{-8}$.

Our algorithm is based on Bandt's algorithm to study the Mandelbrot set for pairs of linear maps \cite{bandt}. The idea is the following: assume that $f\in\mathcal{B}$ has a double zero in $(a,b)\subset (1/2,1)$, and let $P_n\in\mathcal{B}_n$ be the initial part of $f$ up to exponent $n$. Then, letting $r$ be the double root,

\begin{eqnarray}
|P_n(r)| & = & |P_n(r)-f(r)| \le \sum_{i=n+1}^\infty r^i < \frac{b^{n+1}}{1-b}; \nonumber\\
|P'_n(r)| & = & |P'_n(r)-f'(r)| \le \sum_{i=n+1}^\infty i r^{i-1} <  \frac{ (1+n(1-b)) b^n }{ (1-b)^2}.\nonumber
\end{eqnarray}

Hence, using that $a<r<b$ and the intermediate value theorem,

\begin{eqnarray}
|P_n(b)| & \le & |P_n(r)| + \|P'_n\|_{L^{\infty}(r,b)} (b-r) < \frac{b^{n+1}}{1-b} + \frac{ b-a}{ (1-b)^2}; \nonumber \\
|P'_n(b)| & \le & |P'_n(r)| + \|P''_n\|_{L^{\infty}(r,b)} (b-r) < \frac{ (1+n(1-b)) b^n }{ (1-b)^2} + \frac{ b-a}{ 2 (1-b)^3}. \label{eq-boundP}
\end{eqnarray}

Conditions (\ref{eq-boundP}) are easily checkable. If for some $n$ and $(a,b)$, at least one of them fails for all $P\in\mathcal{B}_n$, then $(a,b)\cap \Omega_2 = \varnothing$.

We discuss some of the details of our implementation. Full C++ code is available upon request.

We used floating point arithmetic in order to get practical performance (in principle, one could use exact rational arithmetic instead). To keep the algorithm rigorous we computed the theoretical floating-point error of our calculations and added a corresponding error term to the inequalities (\ref{eq-boundP}). In order to reduce the number of arithmetic operations and avoid very small numbers and the consequent loss of precision we multiplied both sides of (\ref{eq-boundP}) by $b^{-n}$. Thus the numerical checks that we used are
\begin{eqnarray}
b^{-n}|P_n(b)| & < & \frac{b}{1-b} + \frac{(b-a)b^{-n}}{ (1-b)^2} + \eta; \nonumber \\
b^{-n}|P'_n(b)| & < & \frac{ 1+n(1-b) }{ (1-b)^2} + \frac{ (b-a)b^{-n}}{ 2 (1-b)^3} + \eta \label{eq-numboundP}.
\end{eqnarray}

Here $\eta$ is the error term. In practice taking $\eta=10^{-14}$ suffices. This was calculated based on the IEEE floating-point standard. We used the MinGW compiler on a Windows XP platform.

In order to make the algorithm efficient we exploited the tree structure of $\mathcal{B}_n$. The basic routine takes as arguments an interval $(a,b)$, a depth $d$, and a polyomial $P\in\mathcal{B}_n$. The routine returns a boolean value, which we denote by $C((a,b),d,P)$. This value indicates whether the inequalities (\ref{eq-numboundP}) are verified for at least one polynomial of degree at most $n+d$ with initial part $P$. Hence if $C((a,b),d,1)$ returns \textit{false} for some $d$ we must have $(a,b)\cap \Omega_2=\varnothing$. The structure of the routine is as follows:

\begin{enumerate}
\item[(i)] Check (\ref{eq-numboundP}) for $P$ and the interval $(a,b)$. If any of the inequalities fails to hold, return \textit{false} and exit.
\item[(ii)] If $d=0$, then return \textit{true} and exit.
\item[(iii)] For $i=-1, 0$ and $1$, run $C((a,b),d-1, P(x)+i x^{n+1})$. If any of these returns \textit{true}, then return \textit{true} and exit. \label{it-algo}
\item[(iv)] Return \textit{false}.
\end{enumerate}

Assuming that $C((a,b),d,1)$ returns $\textit{true}$, this routine easily produces a polynomial of degree $d$ for which (\ref{eq-numboundP}) holds (this is done by keeping track of which $i$ produces a \textit{true} in Step \ref{it-algo} of the routine). This polynomial can in turn be used to show that there is a double root near $b$ using the results of Section \ref{sec:exist}.

\begin{figure}
\centering
\includegraphics[width=0.9\textwidth]{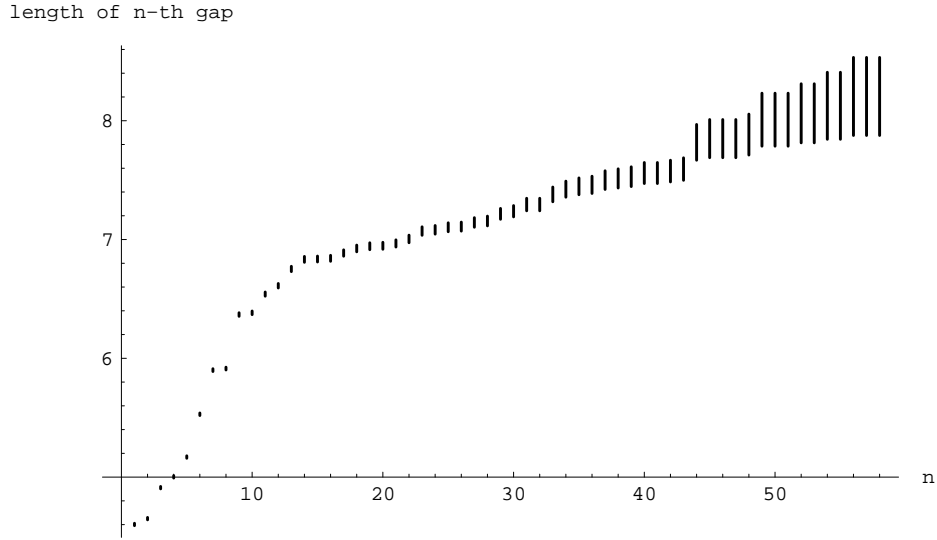}
\caption{This figure shows the distribution of the largest gaps in the set $\Omega_2$. The negative of the base $10$ logarithms of the length of the gaps are contained in the vertical segments (this is a rigorous consequence of the algorithm applied); heuristically, we believe that the actual lengths are closer to the upper end of the intervals.  }
\label{fig-gaps}
\end{figure}

Note that in order to use a large depth one needs to run the algorithm on very small intervals, due to the presence of the term $(b-a)b^{-n}$ in the right-hand side of (\ref{eq-numboundP}). Even on a standard desktop PC, it took less than 3 hours to scan the interval $(0.66847,0.66936)$ for gaps using a grid of $10^7$ subintervals and running the main procedure on each. Figure \ref{fig-gaps} summarizes our findings. We plotted the logarithms of the lengths of the 60 gaps we found, as well as the lengths of the complementary intervals. Although we do not know how to prove it, we believe that at least the few largest of those pieces which the algorithm did not rule out contain intervals of the set $\Omega_2$ of approximately the same length.

We finish this section explaining how we verified the assertions of
Theorem~\ref{th-main}.

\medskip

{\em Proof of Theorem~\ref{th-main} }. We used the algorithm to show that $(0.5,0.6684754)$ and $I_j\, (1\le j\le 5)$ are contained in $(0,1)\backslash\Omega_2$. We did not need to go deeper than 50 iterations for this. For example, we subdivided the interval $I_1$ into $110$ intervals $(a_i,b_i)$ of length $10^{-6}$ and ran $C((a_i,b_i),40,1)$ for each $i$. The routine returned a \textit{false} value for all, which implies that $I_1\in (0,1)\backslash\Omega_2$.

It follows from Example~\ref{ex-double} that $\alpha_2<0.6684757$. To check that the intervals $I_j$ lie in different connected components of $(0,1)\backslash\Omega_2$ we used the algorithm to produce suitable polynomials, and then followed the scheme of Example \ref{ex-double} to show the existence of points $\theta_j (1\le j\le 4)$, $\theta_j \in (\max I_j,\min I_{j+1})\cap \Omega_2$.

The following table summarizes our findings. The symbols on the right represent the coefficients of the polynomials $P_j$; $\o$ corresponds to the coefficient $-1$. The numbers $\theta'_j$ are within $10^{-7}$ of a double root; more precisely, for each $j$ there exist $f_j\in\mathcal{B}$ with initial part $P_j$, and $\theta_j$ a double root of $f_j$ such that $|\theta_j-\theta'_j|<10^{-7}$. \qed

\medskip

\begin{tabular}{|c|c|c|}
  \hline
  $j$ & $\theta'_j$ & $P_j$   \\
  \hline
  $1$ & $0.668550$ & $1{\o}{\o}{\o}10110110111011111{\o}011111{\o}{\o}1{\o}{\o}0111111111111{\o}{\o}1{\o}{\o}$  \\
  $2$ & $0.668900$ & $1{\o}{\o}{\o}1011011011111{\o}111{\o}{\o}11{\o}{\o}11{\o}111{\o}10{\o}{\o}0{\o}0{\o}{\o}{\o}{\o}0{\o}0{\o}{\o}{\o}$  \\
  $3$ & $0.669310$ & $1{\o}{\o}{\o}101101110011101011110111111{\o}000{\o}{\o}{\o}{\o}{\o}1{\o}{\o}{\o}{\o}100011$  \\
  $4$ & $0.669336$ & $1{\o}{\o}{\o}1011011100111011001111111{\o}0110101111{\o}1111111{\o}{\o}0$  \\
  \hline
\end{tabular}

\medskip

{\em Proof of Corollary~\ref{cor-ant}.} The double zero $y$ obtained in
Example~\ref{ex-double} is a root of $f\in \Bk$ which has infinitely many 0's
among the coefficients, by Remark~\ref{rem-cor}. Replacing
$a_j =0$ by $+1$ for $j$ large yields a function in $\Bk$ with  a complex
zero close to $y$. On the other hand, replacing $a_j =0$ by $-1$ for
$j$ large yields a function in $\Bk$ with  two real zeros close to $y$. Thus,
$y \in \clos(\Mk\setminus \R)$ and
$(y,y) \in \clos(\Nk\setminus \Diag(\R))$. In view of Lemma~\ref{lem-anten},
this implies both ``tips of antennas'' belong to $(\alpha_2,y)$, a
and the claim follows from Theorem~\ref{th-main}. \qed


\section{Variants and generalizations}

The algorithm described in the previous section adapts without difficulty to
more general settings; we consider some examples below.

\subsection{The set of triple roots}

Here we restrict ourselves to the family $\mathcal{B}$, but consider higher order roots. Geometrically, the set of roots of multiplicity $n$ corresponds to the connectedness locus of self-affine sets associated to Jordan blocks of order $n$; recall Proposition \ref{prop-connect}. We denote this set by $\Omega_n$.

The computer algorithm extends in a straightforward way to higher multiplicity roots. Indeed, it is enough to replace (\ref{eq-boundP}) by the set of tests
\[
|P^{(i)}(b)| < H_n^{(i)}(b) + \frac{b-a}{(i+1) (1-b)^{i+2}}, \quad 0\le i< n,
\]
where $H_n(x) = x^n/(1-x)$.
The algorithm yields intervals in $(0,1)\backslash\Omega_n$; in particular, it gives lower bounds on $\alpha_n=\min\Omega_n$. We remark, however, that in practice the program becomes very slow for $n\ge 4$.  For $n=3$, we have the following result.

\begin{prop}
{\bf (i)}
$\alpha_3=\min\Omega_3 > 0.743$,

{\bf (ii)} $(0.746,0.7465) \subset (0,1)\backslash\Omega_3$,

{\bf (iii)} $(2^{-1/3},1)\subset\Omega_3$.
\end{prop}

{\em Remark on proof:} (i) and (ii) are direct applications of the algorithm,
while (iii) follows from Lemma \ref{lem-folk}.
(Note that $2^{-1/3}\approx 0.7937$). \qed

\begin{remark}
Numerical experimentation suggests that $\alpha_3\in (0.743,0.744)$ and, in particular, $\Omega_3$ is disconnected. However, the techniques of section \ref{sec:exist} do not seem to apply, so this remains a conjecture.
\end{remark}

It is interesting to compare our results with those of \cite{BBBP},
where multiple roots of the following family were considered:
$
\widetilde{\mathcal{B}} = \left\{ 1 + \sum_{i=1}^\infty a_i x^i : a_i \in [-1,1] \right\}.
$
Let $\beta_n$ be the smallest root of multiplicity (at least) $n$ of some $f\in \widetilde{\mathcal{B}}$. The values of $\beta_n$ for $n\le 27$ were computed in \cite{BBBP}; in particular, $\beta_2\approx 0.64914$ and
$\beta_3\approx 0.72788$. Observe that $\alpha_2-\beta_2 > 0.01934$ and
$\alpha_3-\beta_3>0.0151$. Thus going from a continuous to a discrete set of
coefficients does have a substantial impact in the set of multiple roots.

\subsection{The set of double zeros with coefficients $0,\pm 1,\pm 2$}

We can generalize the set $\Omega_2$ in another direction by enlarging the set of allowed coefficients. For concreteness, we will work with the coefficient set $\{0,\pm 1,\pm 2\}$; specifically, let
\[
\mathcal{B}' = \Bigl\{ 1 + \sum_{i=1}^\infty a_i x^i : a_i \in \{0,-1,1,-2,2\}
\Bigr\}.
\]

Denote by $\Om_2'$ the set of double zeros in $(0,1)$ of elements  of
$\mathcal{B}'$. It turns out that a specific power series plays a very special
role in the study of this set. Let
\[
Q(x) = 1 - 2x - 2x^2 + \sum_{i=3}^\infty 2 x^i = 1 - 2x - 2x^2 +
\frac{2x^3}{1-x}\,.
\]
Note that $1/2$ is a double root of $Q(x)$. In fact, more is true: $Q(x)$ is a so-called $(*)$-function for the class $\mathcal{B}'$ on the interval $(0,x_0)$ for all $x_0<1/2$; see \cite{solnotes} for the relevant definitions and proofs. A consequence of this is that $\min\Omega'_2$ is precisely $1/2$.

Here we prove that $1/2$ is actually an isolated point of $\Omega'_2$; this is due to the special form of the function $Q(x)$, and in particular the fact that all but finitely many coefficients are $+2$. More precisely, we have the following result:

\begin{prop} $\min(\Omega'_2\backslash\{1/2\}) \in (0.5436,0.5438)$. \label{prop-isolated}
\end{prop}

Before proving the proposition we remark that the set-up of Section
\ref{sec:exist} works here with minor modifications. In this context, we say
that $(P,n,a,b)$ is {\em good} if $P\in\mathcal{B}'_n$
(the family of polynomials in $\mathcal{B}'$ of degree at most $n$), $0.5<a<b<1$,
\[
P(a) > 2 a^{n+1}/(1-a),\quad P(b) > 2 b^{n+1}/(1-b),
\]
$P(x)>0$ for all $x\in[a,b]$, and
\[
\exists\,x\in (a,b):\ P(x) < 2 x^{n+1}/(1-x).
\]
The proofs of Section \ref{sec:exist} apply almost verbatim. In particular, if $(P,n,a,b)$ is good then there exists a sequence $i_1<i_2<\ldots$ such that $f(x)=P(x)-\sum_{j=1}^\infty 2 x^{i_j}$ has a double root in $(a,b)$.

\medskip

{\em Proof of Proposition \ref{prop-isolated}.} We will use the following result, which follows from a modification of the proof of Theorem 2 in \cite{BBBP}
(see \cite{shm} for a complete proof): if $f$ is a power series with coefficients in $[-1,1]$, and $\alpha_1,\ldots, \alpha_k$ are complex roots of $f$ in the unit disk, counted with multiplicity, then
\begin{equation} \label{eq-triplerootbound}
|\alpha_1\,\ldots\,\alpha_k| \ge \left(1+\frac{1}{k}\right)^{-k/2} (k+1)^{-1/2}.
\end{equation}
We will use this result with $k=3$. A standard application of the algorithm shows that $f\in\mathcal{B}'$ may have roots in $(0.5,0.51)$ only if it starts
with $1-2x-2x^2$. Then it follows from the definition of $Q$ that $f(x) < Q(x)$
for all $x>0$. Since $f(0)=1$ and $f(1/2)<Q(1/2)=0$, $f$ has a root in the interval $(0,1/2)$. Suppose that $\alpha$ is a double root of $f$ in $(1/2,1)$. We obtain from (\ref{eq-triplerootbound}) that
\[
\frac{1}{2} \alpha^2 > \left(1+\frac{1}{3}\right)^{-3/2} (3+1)^{-1/2} = \frac{3\sqrt{3}}{16}\,,
\]
whence $\alpha > (6\sqrt{3}/16)^{1/2} > 0.8$. We conclude that $f$ cannot have
double roots in the interval $(0.5,0.51)$, whence $(0.5,0.51)\subset
(0,1)\backslash\Omega'_2$.

A standard application of the algorithm shows also that
$(0.51,0.5436)\subset(0,1)\backslash\Omega'_2$.

Finally, let $P$ be the polynomial of degree $26$ with coefficients
\[
(1, -2, -1, 1, 1, 1, 2, 1, 1, 2, 1, 1, 2, 1, 1, 2, 1, 2, 1, 1, 2, -2, -2, -1, -2, 2, -1)
\]
We checked with {\em Mathematica} that $(P,27,0.5436,0.5438)$ is good.
This implies that $\min\Omega'_2 \in (0.5436,0.5438)$,
completing the proof. \qed

\begin{remark}
The set $\Omega'_2\backslash\{1/2\}$ seems to be connected (i.e.\ no ``gaps''
appear when running the program), but we do not have a proof of this.
Still, if what the numerical experimentation suggests holds true, then the
sets $\Omega_2$ and $\Omega'_2$ have strikingly different
topological structure. \end{remark}


\section{Remaining proofs}

{\em Proof of Lemma~\ref{lem-folk}.} Let $E$ be the attractor of the IFS,
that is, $E = T_1E + (T_2E + \bob)$.  Let $\|\cdot\|$ be the norm in
$\R^d$ such that $\|T_ix\| \le r\|x\|,\ i=1,2,$
for some $r \in (0,1)$ and all $x\in \R^d$.
Denote by $F_\eps$ the $\eps$-neighborhood of a
set $F\subset \R^d$ in this norm. Then $(T_iF)_\eps \supset
T_i F_{\eps/r}$.

Suppose that $E$ is
disconnected. Then $T_1E \cap (T_2E+\bob) = \es$ by
Proposition~\ref{prop-hata}, and since these are compact sets,
we can find $\eps>0$ such that $(T_1E)_{\eps} \cap (T_2E + \bob)_\eps  = \es$.
Then $E_\eps = (T_1E)_\eps \cup (T_2E + \bob)_\eps$ is a disjoint union,
so
\begin{eqnarray*}
\Leb^d(E_\eps) & = & \Leb^d((T_1E)_\eps) + \Leb^d((T_2E)_\eps) \\
& \ge &  \Leb^d(T_1 E_{\eps/r}) + \Leb^d(T_2 E_{\eps/r}) \\
& = & (|\det(T_1)|+|\det(T_2)|) \cdot
\Leb^d(E_{\eps/r}) \ge \Leb^d(E_{\eps/r}).
\end{eqnarray*}
This is a contradiction, since $E_{\eps/r}\setminus E_\eps$ has
positive Lebesgue measure. \qed

\medskip

{\em Proof of Proposition~\ref{prop-connect}.}
If $E= E(T,\bob)$ is connected, then $TE \cap (TE + \bob)\ne \es$.
In view of (\ref{serep}), we obtain that there exist
$\{0,1\}$ sequences $\{a_n\}_0^\infty$ and $\{a_n'\}_0^\infty$
such that $a_0=1, a'_0=0$, and
$ \sum_{n=0}^\infty (a_n - a_n') T^n \bob = \bo$.
Denoting $f(x) =  \sum_{n=0}^\infty (a_n - a_n') x^n$ we get
$f \in \Bk$ and $f(T) \bob = \bo$. Now let us write
$\bob = \sum_{j=1}^m c_j \bbe_j$ where $\bbe_j \in \Ker(T-\lam_j I)^{k_j}$.
Since $\bob$ is a cyclic vector for $T$, we have $c_j \ne 0$ and $\bbe_j
\not\in \Ker(T-\lam_j I)^{k_j-1}$ for $j\le m$.
Then $f(T)\bob= \bo$ implies that $f(T) \bbe_j = \bo$ for all $j\le m$.
We have
$$
f(T) \bbe_j = f(\lam_j) \bbe_j + f'(\lam_j) (T-\lam_j I) \bbe_j +\ldots
+ f^{(k_j-1)}(\lam_j) (T - \lam_j I)^{k_j-1} \bbe_j.
$$
Since the vectors $\{(T-\lam_j I)^\ell \bbe_j:\ \ell = 1,\ldots, k_j -1\}$
are linearly independent, (\ref{eq-zero}) follows.

Conversely, if $f(x)=1+\sum_{n=1}^\infty b_n x^n\in \Bk$
satisfies (\ref{eq-zero}), then  $f(T)\bob = \bo$ for all $\bob$.
Writing $b_n = a_n - a_n'$ for some $a_n, a_n' \in \{0,1\}$, we obtain that
$TE \cap (TE + \bob)\ne \es$,
hence $E$ is connected by Proposition~\ref{prop-hata}.
\qed

\medskip

{\em Proof of Lemma~\ref{lem-anten}.} (ii) Suppose that $\lam\in (-1,1)$ and
$(\lam,\lam)$ is
such that there exists a sequence $(\gam_n,\lam_n)\in \Nk$,
with $\gam_n< \lam_n$, converging to $(\lam,\lam)$. Then there are
power series $f_n \in \Nk$ such that $f_n(\lam_n) = f(\gam_n) =0$.
By compactness, passing to a subsequence,
we can assume that $f_n \to f\in \Bk$ coefficientwise. Then $f^{(k)}_n(\lam)\to
f^{(k)}(\lam)$ for all $k\ge 0$. By Lemma~\ref{lem-folk}, $|\lam| < 2^{-1/2}$.
Let $C_k = \max \{|f^{(k)}(x)|:\ f\in \Bk,\ |x|\le 2^{-1/2}\}$, which is
finite (and easy to compute explicitly). Then we have for $n$ sufficiently
 large:
$$
|f_n(\lam)| = |f_n(\lam) - f_n(\lam_n)| \le C_1|\lam-\lam_n| \to 0,\ \ n\to
\infty.
$$
Next, there exists $t_n \in (\gam_n,\lam_n)$ such that $f'_n(t_n)=0$ and
we have for $n$ sufficiently  large:
$$
|f'_n(\lam)| = |f'_n(\lam) - f'_n(t_n)| \le C_2|\lam-t_n| \to 0,\ \ n\to
\infty.
$$
It follows that $f(\lam) = f'(\lam)=0$, as desired. \qed

\medskip

{\em Proof of Lemma~\ref{lem-vsp1}.}
In the following calculations $\eta$ will be a fixed positive number,
to be determined later. We use the $\ell^\infty$ norm on $\R^2$.

Given a word $u\in \{-1,0,1\}^5$ let
\begin{eqnarray}
T_u & = & (T+u_1 \bob)\circ\cdots\circ (T+u_5 \bob),\nonumber\\
T'_u & = & (T'+u_1 \bob)\circ\cdots\circ (T'+u_5 \bob).\nonumber
\end{eqnarray}
Observe that
\begin{eqnarray}
T_u \bx & = & T^5 \bx + \sum_{i=1}^5 u_i T^{i-1} \bob\nonumber,\\
(T'_u)^{-1}\bx & = & (T')^{-5}\bx -
(T')^{-5}\left(\sum_{i=1}^5 u_i (T')^{i-1}\bob\right).\nonumber
\end{eqnarray}
Hence
\begin{eqnarray}
(T'_u)^{-1}T_u  \bx & = & \left((T')^{-5}T^{5}\right) \bx + (T')^{-5}
\sum_{i=1}^5 u_i T^{i-1} \bob -(T')^{-5} \sum_{i=1}^5 u_i (T')^{i-1}\bob
\nonumber\\
& =: & S\bx + \bd,\label{eq-normest1}
\end{eqnarray}
where
\[
S = (T')^{-5} T^5,\ \ \ \ \bd = (T')^{-5}\sum_{i=1}^5 (T^{i-1} - (T')^{i-1})
\bob.
\]
Observe that (for $\eta < 10^{-2}$)
\begin{eqnarray}
\|T^k - (T')^k\| & \le & \sum_{j=1}^k \|T^j (T')^{k-j} - T^{j-1} (T')^{k-j+1}\| \nonumber \\
& = & \sum_{j=1}^k \| T^{j-1} (T-T') (T')^{k-j-1}\| \nonumber\\
& \le & k \eta \max(\|T\|,\|T'\|)^{k-1} <  k \,(1.5)^{k-1}\,\eta. \label{dub}
\end{eqnarray}
Hence if $R=(T')^5-T^5$ then $\|R\|\le 5\,(1.5)^4\,\eta < 26 \eta$,
and we can estimate
\begin{eqnarray}
\|(T')^{-5}\| & = & \|\left(T^{5}(I-T^{-5} R)\right)^{-1}\| \le
\|(I-T^{-5}R)^{-1}\| \|T^{-5}\| \nonumber\\
& \le & \|T^{-5}\| \sum_{j=0}^\infty (\|T^{-5}\| \|R\|)^j =
\frac{\|T^{-5}\|}{1-\|T^{-5}\|\|R\|} \nonumber\\
& < & \frac{33.6}{1-33.6\times 26\,\eta} < 34,
\end{eqnarray}
as long as $\eta<10^{-5}$. Therefore,
\[
\|S-I\| = \|(T')^{-5}R\| \le \|(T')^{-5}\| \|R\| < 34 \times 26\,\eta =
884 \,\eta.
\]
We have by (\ref{dub}),
\begin{equation} \label{dub-est}
\|\bd\| \le \|(T')^{-5}\| \sum_{j=0}^4 j\,(1.5)^{j-1}\eta \|\bob\| < 34 \times
24.25\,\eta < 825\, \eta.
\end{equation}
Note that
\[
\max\{\|\bx\| :\bx\in V \} = 0.95 \max(\|\textbf{p}\|,\|\textbf{q}\|) =
0.95 \times 4.9 < 4.7.
\]
It follows from the previous estimates that for $x\in V$,
\[
\|S \bx - \bx\| \le \|S-I\| \|\bx\| < 884\times 4.7\, \eta < 4155\,\eta.
\]
In particular, this implies that $S V \subset V_\delta$, where
$\delta = 4155\, \eta$, and $V_\delta$ denotes the $\delta$
neighborhood of $V$. Recalling (\ref{eq-normest1}) we get that
\begin{equation} \label{eq-normestsubset1}
(T'_u)^{-1}T_u  V \subset (V_\delta)_{\|\bd\|} = V_{\delta+\|\bd\|}
\subset V_{\delta'},
\end{equation}
where $\delta' = 5 \times 10^3\,\eta$, since $\|\bd\|<825\, \eta$ by
(\ref{dub-est}).

Let $W(r)=\{(x,y):|x|+|y|< r\}$ and let $M$ be the matrix with the
columns $\mathbf{p},\mathbf{q}$. Note that $U=M W(1)$ and $V=M W(0.95)$.
Note also that
\[
\dist(W(0.95),\mathbb{R}^2\backslash W(1)) = 0.05/2 = 2.5 \times 10^{-2},
\]
where $\dist(\cdot,\cdot)$ denotes the distance induced by the
$\ell^\infty$ norm. An easy calculation yields $\|M^{-1}\|=1.2$. Therefore,
\begin{equation}  \label{eq-normestsubset2}
\dist(V,\mathbb{R}^2\backslash U) \ge \frac{\dist(W(0.95),\mathbb{R}^2\backslash W(1))}{\|M^{-1}\|} > 2 \times 10^{-2}.
\end{equation}
It follows from (\ref{eq-normestsubset1}) and (\ref{eq-normestsubset2})
that if $\eta$ is so small that $\delta'\le 2 \times 10^{-2}$, then
\begin{equation}  \label{eq-normestsubset3}
(T'_u)^{-1} T_u V \subset U.
\end{equation}
Since $\delta' = 5 \times 10^3 \eta$, this will be the case for
$\eta = 4 \times 10^{-6}$. From now on, we will fix this value of $\eta$
(since $\eta<10^{-5}$, the previous calculations apply).

Let $A$ be the ``multiplication by $0.95$'' map, so that $V=AU$.
By applying $A$ to both sides of (\ref{eq-parallcover}) we get
\[
V \subset A \left(\bigcup_{u\in\{-1,0,1\}^5}  T_u(V)\right) = \bigcup_{u\in\{-1,0,1\}^5} T_u A V.
\]
However, we deduce from (\ref{eq-normestsubset3}) that $T_u V \subset
T'_u A^{-1} V$, or $T'_u V \supset T_u A V$. Combining this with the last
displayed formula, we conclude that (\ref{eq-claimcovering})  holds. \qed


\bibliographystyle{amsplain}

\end{document}